\documentclass[11pt]{article}

\usepackage{amsxtra,amssymb,amsthm,amsmath,latexsym}
\usepackage{graphicx}
\usepackage{epsfig}
\usepackage{afterpage}
\newtheorem{thm}{Theorem}[]

\newtheorem{lem}[]{Lemma}
\newtheorem{dfn}[]{Definition}

 \newcommand{\thmref}[1]{Theorem~\ref{#1}}
 \newcommand{\lemref}[1]{Lemma~\ref{#1}}

\newcommand{\R}{{\mathbb R}}

\newcommand{\dl}{{\delta}}
\newcommand{\bee}{\begin{equation*}}
\newcommand{\eee}{\end{equation*}}
\newcommand{\be}{\begin{equation}}
\newcommand{\ee}{\end{equation}}
\newcommand{\pn}{\par\noindent}

\title{Stability of solutions to abstract differential equations}
\author{A.G. Ramm \\
\small Department of Mathematics\\[-0.8ex]
\small Kansas State University, Manhattan, KS 66506-2602, USA\\
\small \texttt{ramm@math.ksu.edu}}

\date{}
\begin{document}

% typeset front matter
\maketitle

\begin{abstract}
A sufficient condition for asymptotic stability of the zero solution
to an abstract nonlinear evolution problem is given. The governing
equation is $\dot{u}=A(t)u+F(t,u),$ where $A(t)$ is a bounded linear
operator in Hilbert space $H$ and $F(t,u)$ is a nonlinear operator,
$\|F(t,u)\|\leq c_0\|u\|^{1+p}$, $p=const >0$, $c_0=const>0$. It is not 
assumed
that the spectrum $\sigma:=\sigma(A(t))$ of $A(t)$ lies in the fixed 
halfplane Re$z\leq -\kappa$, where  $\kappa>0$ does not depend on $t$. As 
$t\to \infty$
the spectrum of $A(t)$ is allowed to tend to the imaginary axis.
\end{abstract}
\pn{\\ {\em MSC:} 34G20; 447J05; 47J35  \\

\noindent\textbf{Key words:} dynamical systems; stability;  asymptotic stability }\\

\section{Introduction}
Let $H$ be a Hilbert space. Consider the problem
 \be\label{e1}
\dot{u}=A(t)u+F(t,u),\quad t\geq 0, \ee \be\label{e2} u(0)=u_0,\ee
where  $\dot{u}=\frac{du}{dt}$ is the strong derivative, $A(t)$ is a
linear closed densely defined in $H$ operator with the domain
$D(A)$, independent of $t$, $u_0\in D(A)$. 
We assume that $F(t,u)$ is a nonlinear
mapping, locally Lipschitz with respect to $u$, and satisfying the 
following inequality
 \be\label{e3}
\|F(t,u)\|\leq c_0\|u\|^{1+p},\quad p>0,\ c_0>0, \ee 
where $p$ and $c_0$
are constants. We also assume that ~\be\label{e4} \text{Re}(Au,u)\leq
-\gamma(t)\|u\|^2,\quad \forall u\in D(A), \ee where \be\label{e5}
\gamma(t)>0, \quad \lim_{t\to \infty}\gamma(t)=0, \ee \be\label{e6}
\gamma(t)=\frac{b_1}{(b_0+t)^d},\quad d=const\in(0,1], \ee $b_0$ and
$b_1$ are positive constants. Assumptions \eqref{e5} are satisfied by
the function \eqref{e6}. However, our method can be applied 
to many other $\gamma(t)$ satisfying assumptions \eqref{e5}.

\begin{dfn}\label{def1}
The zero solution to equation \eqref{e1} is called Lyapunov stable
if for any $\epsilon>0$, sufficiently small, there exists a
$\dl=\dl(\epsilon)>0$, such that if $\|u_0\|<\dl$, then the solution to
problem \eqref{e1} exists on $[0,\infty)$ and $\|u(t)\|\leq
\epsilon$. If, in addition, \be\label{e7} \lim_{t\to
\infty}\|u(t)\|=0, \ee then the zero solution is asymptotically
stable.
\end{dfn}

Basic results on the Lyapunov stability of the solutions to
\eqref{e1} one finds in \cite{Ba}-\cite{DK},  and in many other
books and papers. In \cite{DK} these results are established under
the assumption that the operator $A(t)$ is bounded, $D(A)=H$, and
$A(t)$ has property $B(\nu,N)$. This means (\cite{DK}, p.178) that
every solution to the equation \be\label{e8} \dot{u}=A(t)u \ee
satisfies the estimate \be\label{e9} \|u(t)\|\leq
Ne^{-\nu(t-s)}\|u(s)\|,\quad t\geq s\geq 0, \ee where $N>0$ and
$\nu>0$ are some constants.  The quantity \be\label{e10}
\kappa:=\overline{\lim_{t\to \infty}}\frac{\ln \|u(t)\|}{t} \ee is
called the exponent of growth of $u(t)$. If $\Sigma$ is the set of
$\kappa$ for all solutions to \eqref{e8}, then \be\label{e11}
\kappa_s:=\sup_{\kappa\in \Sigma}\kappa \ee is called senior
exponent of growth of solutions to \eqref{e8}. The general exponent
$\kappa_g$ is defined as \be\label{e12} \kappa_g:=\inf \rho, \ee
where $\rho$ is the exponent in the inequality \be\label{e13}
\|u(t)\|\leq Ne^{\rho(t-s)}\|u(s)\|,\quad t\geq s\geq 0. \ee One has
\be\label{e14} \kappa_s\leq \kappa_g, \ee and the case
$\kappa_s<\kappa_g$ can occur (the Perron's example, see \cite{DK},
p.177). If $\kappa_g<0$ then the zero solution to \eqref{e8} is
Lyapunov asymptotically stable. If $A(t)=A$ does not depend on $t$
and $A$ is a bounded linear operator, then $\kappa_g<0$  if and only
if the spectrum of $A$, denoted $\sigma(A)$, lies in the halfplane
Re$z\leq \kappa_g<0$. In this case \be\label{e15} \|e^{At}\|\leq
N_0e^{\kappa_gt}, \ee and if $\|F(t,u)\|\leq q\|u\|,$ $t\geq 0$,
$\|u\|<\rho$, and $q<\frac{\kappa_g}{N_0}$, then equation \eqref{e1}
has negative general exponent also, so the zero solution to equation
\eqref{e1} is Lyapunov asymptotically stable (\cite{DK}, p.403).

If
$A=A(t)$, and for any solution to \eqref{e8} estimate \eqref{e9}
holds with $\nu>0$, and if \eqref{e3} holds, then for any solution
to \eqref{e1} with $\|u_0\|\leq \dl$ and $\dl>0$ sufficiently small,
estimate \eqref{e9} holds with a different $N=N_1$ and $\nu=\nu_1$,
$0<\nu_1\leq \nu$ (see \cite{DK}, p.414). This means that the zero
solution to \eqref{e1} is asymptotically stable under the above
assumptions.

The basic new result of our work, Theorem 1 in Section 2,  generalizes the
above results to the case when the assumption $\kappa_g<0$ is not valid.
We allow the spectrum $\sigma(A(t))$ to approach imaginary axis
as $t\to \infty$. This is the principally new generalization of the
classical Lyapunov-Krein theory.  If
$\sqcap$ is the complex plane and $l$ is the imaginary axis, then
we assume that $\sigma(A(t))\subset \sqcap$ for every $t\geq 0$, but we 
allow 
$\lim_{t\to\infty}d(\sigma(A(t)),l)=0,$ where $d(\sigma, l)$ is the 
distance between two sets $\sigma$ and $l$.
The new stability result is formulated
in \thmref{thm1}. In \lemref{lem1} an auxiliary result is formulated.
A proof of \lemref{lem1} differs in details from the one in 
\cite{R558}.
In Section 2 \thmref{thm1} and \lemref{lem1} are formulated. In
Section 3 proofs are given. In Section 4 examples of applications
of our method are given.

\section{Formulation of the results}
\begin{lem}\label{lem1}
Let the inequality \be\label{e16} \dot{g}(t)\leq
-\gamma(t)g(t)+a(t)g^{1+p}(t)+\beta(t), \ee hold for
$t\in[0,T)$,where $g(t)\geq 0$ has finite derivative from the right at
every
point $t$ at which $g(t)$ is defined, $\gamma(t)\geq 0$, $a(t)\geq
0$ and $\beta(t)\geq 0$ are continuous on $\R_+:=[0,\infty)$
functions, and $p=const>0$. Assume that there exists a $\mu(t)\in
C^1[0,\infty)$,
$\mu(t)>0,$ $\dot{\mu}(t)\geq 0$, such that
\be\label{e17}
a(t)[\mu(t)]^{-1-p}+\beta(t)\leq
\mu^{-1}(t)[\gamma(t)-\dot{\mu}(t)\mu^{-1}(t)],\ t\geq 0, \ee
\be\label{e18} \mu(0)g(0)<1. \ee Then $g(t)$ exists for all
$t\in [0,\infty)$
and \be\label{e19} 0\leq g(t)<\mu^{-1}(t),\quad \forall t\geq 0. \ee
\end{lem}

\begin{thm}\label{thm1}
Assume that conditions  \eqref{e1}-\eqref{e6} hold
and $b_1>0$ is sufficiently large. Then the zero solution
to
\eqref{e1} is asymptotically stable for  any fixed initial data $u(0)$.
\end{thm}
\section{Proofs}
\noindent {\it Proof of \lemref{lem1}}. Let
$v(t):=g(t)e^{\int_0^t\gamma(s)ds}:=g(t)q(t).$ Then \eqref{e16}
yields
\be\label{e20} \dot{v}(t)\leq
q(t)a(t)q^{-(1+p)}(t)v^{1+p}(t)+q(t)\beta(t),\quad v(0)=g(0),\quad
t>0. \ee
We do not assume a priori  that $v(t)$ is defined for
all $t\geq 0$. This conclusion will follow from our proof.
Denote $\eta(t):=q(t)\mu^{-1}(t)$, $\eta(0)=\mu^{-1}(0)>g(0)$. Using
\eqref{e18} and \eqref{e20}, one gets \be\label{e21} \dot{v}(0)\leq
a(0)v^{1+p}(0)+\beta(0)\leq
\mu^{-1}(0)[\gamma(0)-\dot{\mu}(0)\mu^{-1}(0)]=\dot{\eta}(0). \ee
Since $v(0)=g(0)<\eta(0)=\mu^{-1}(0)$ by \eqref{e18}, and $\dot{v}(0)\leq 
\dot{\eta}(0)$, it
follows that \be\label{e22} v(t)<\eta(t),\quad 0\leq t<\tau, \ee
where $\tau>0$ is the right end of the maximal interval
on which $ v(t)<\eta(t)$, i.e., $\tau=\sup_{\{t\ :
\ v(t)<\eta(t)\}}t$. Let us prove that
$\tau=\infty$. Note that if \eqref{e22} holds, then \be\label{e23}
\dot{v}(t)\leq \dot{\eta}(t),\quad 0\leq t<\tau. \ee Indeed, using
\eqref{e17} and \eqref{e20} one obtains \be\label{e24}
\dot{v}(t)=q(t)(\dot{g}+\gamma g)\leq
q(t)\mu^{-1}(t)[\gamma(t)-\dot{\mu}(t)\mu^{-1}(t)]=\dot{\eta}(t),
\ee as claimed. If $\tau<\infty$, then \eqref{e22} and \eqref{e23}
imply \be\label{e25} v(\tau-0)-v(0)\leq \eta(\tau-0)-\eta(0). \ee
Since $\eta(t)\in C^1[0,\infty)$ by definition, inequality
\eqref{e25} implies that $v(\tau-0)<\infty$ and, since
$v(0)=g(0)<\mu^{-1}(0)=\eta(0),$ so that $v(0)<\eta(0),$ one gets
\be\label{e26}v(\tau-0)<\eta(\tau-0)<\infty.\ee Inequality
\eqref{e26} implies that $\tau=\infty$, because
$\tau$ is the maximal interval $[0,\tau)$ of the existence of $v$, and
if $\tau<\infty$ is
the right end of the maximal interval of the existence of $v$ then
$\overline{\lim}_{t\to \tau-0}v(t)=\infty$, which contradicts
\eqref{e26}. Thus, $\tau=\infty$ and, therefore, $T=\infty$.\\
\lemref{lem1} is proved.\hfill $\square$\\

\noindent {\it Proof of \thmref{thm1}.} Let $\|u(t)\|=g(t)$.
Multiply \eqref{e1} by $u(t)$,  take the real part, and get
\be\label{e27} g(t)\dot{g}(t)\leq -\gamma g^2(t)+c_0g^{2+p}(t). \ee
Since $g\geq 0$, inequality \eqref{e27} is equivalent to
\be\label{e28} \dot{g}(t)\leq -\gamma(t)g(t)+c_0g^{1+p}(t). \ee If
$g(t)>0$, then \eqref{e28} is obviously equivalent to \eqref{e27}.
If $g(t)=0$ $\forall t\in \Delta$, where $\Delta \subset \R_+$ is an
open set, then $u(t)=0$ $\forall t\in \Delta$, so $u(t)=0$ $\forall
t\geq 0$ by the uniqueness of the solution to the Cauchy problem for
equation \eqref{e1}. This uniqueness holds due to the assumed local
Lipschitz condition for $F$. If $g(t_0)=0$, but $g(t)\neq 0$ for
$(t_0, t_0+\dl)$ for some $\dl>0$, then one divides \eqref{e27} by
$g(t)$ for $t\in (t_0,t_0+\dl)$, then one passes to the limit $t\to
t_0+0$ and gets \eqref{e28} at $t=t_0$. Let us explain the meaning of
$\dot{g}(t_0)$ at a point where $u(t_0)=0$. The function
$\dot{u}(t)$ is continuous and it is known that
$\frac {d\|u(t)\|}{dt}\leq \|\dot{u}(t)\|.$ We
define $\dot{g}(t_0)=\lim_{s\to +0}\|u(t_0+s)\|s^{-1}. $ This 
limit exists and is equal to $\|\dot{u}(t_0)\|$. Choose
\be\label{e29} \mu(t)=\mu(0)e^{\frac{1}{2}\int_0^t\gamma(s)ds},\quad
\dot{\mu}(t)\mu^{-1}(t)=\gamma(t)/2. \ee
{\bf Remark 1.} {\it Note that $\lim_{t\to \infty}\mu(t)=\infty$  if and 
only 
if $\int_0^\infty \gamma(t)dt=\infty$. If $\lim_{t\to 
\infty}\mu(t)=\infty$, then $\lim_{t\to \infty}\|u(t)\|=0$.
Under the assumption \eqref{e6} one has $\int_0^\infty 
\gamma(t)dt=\infty$, and we use this to derive some results about 
asymptotic stability. If $d>1$ in \eqref{e6}, then 
$\int_0^\infty\gamma(t)dt<\infty$, and our methods can be used 
for a derivation of some results on stability, rather than 
asymptotic stability.}

Condition \eqref{e18} is
satisfied if \be\label{e30} \mu(0)<[g(0)]^{-1}, \ee and we choose
$\mu(0)$ so that this inequality holds. Using \eqref{e29}, one sees
that inequality \eqref{e17} is satisfied if \be\label{e31}
2c_0\mu^{-p}(0)\leq
\gamma(t)e^{\frac{p}{2}\int_0^t\gamma(s)ds}, \quad \forall
t\geq 0. \ee Inequality \eqref{e31} is satisfied if \be\label{e32}
2c_0\mu^{-p}(0)\leq \gamma(0),\ee provided that \be\label{e33}
\gamma(0)\leq \gamma(t)e^{\frac{p}{2}\int_0^t\gamma(s)ds} \qquad 
\forall t\geq 0.\ee 
Let us first
use assumption \eqref{e6} with $d\in (0,1)$: \be\label{e34}
\int_0^t\gamma(s)ds=b_1\frac{(b_0+t)^{1-d}-b_0^{1-d}}{1-d},\quad
0<d<1. \ee In this case $\gamma(0)=b_1b_0^{-d},$ and inequality
\eqref{e33} holds if \be\label{e35} 2d<pb_1b_0^{1-d}. \ee 
Inequality
\eqref{e35} is a sufficient condition for the function on the right
of \eqref{e33} to have non-negative derivative for all $t\geq 0$, i.e.,
to be monotonically growing on $[0,\infty)$, if $\gamma(t)$ is
defined in \eqref{e6}. Conditions \eqref{e32} and \eqref{e35} hold if
\be\label{e36}
2c_0\mu^{-p}(0)\leq b_1b_0^{-d} \quad \text{ and }\quad
2d<pb_1b_0^{1-d}.\ee
For any fixed four parameters $d, c_0, p,$ and $\mu(0)<[g(0)]^{-1}$, 
where
$d\in(0,1)$, $c_0>0$, $p>0$, and $\mu(0)>0$,
inequalities \eqref{e36} can be satisfied by choosing {\it sufficiently
large $b_1>0$}. With the choice of $\mu(t)$, given in \eqref{e29},
and the parameters $\mu(0)$, $b_0$ and $b_1$, chosen as above, one
obtains inequality \eqref{e19}: \be\label{e37} 0\leq
g(t)<\frac{e^{-\frac{b_1}{2(1-d)}[(b_0+t)^{1-d}-b_0^{1-d}]}}{\mu(0)},\qquad
d\in(0,1). \ee Since $g(t)=\|u(t)\|$, inequality \eqref{e37} implies
asymptotic stability of the zero solution to equation \eqref{e1} for
any initial value of $u_0$, that is {\it global} asymptotic stability.
Moreover, \eqref{e37} gives a rate of convergence of $\|u(t)\|$ to
zero as $t\to \infty$.

Consider now the case $d=1$, $\gamma(t)=b_1(b_0+t)^{-1}$, 
\be\label{e38}
\int_0^t\gamma(s)ds=b_1\ln\frac{b_0+t}{b_0},\qquad
e^{\int_0^t\gamma(s)ds}=\left( \frac{b_0+t}{b_0} \right)^{b_1}. \ee In 
this
case the choice of $\mu(t)$ in \eqref{e29} yields \be\label{e39}
\mu(t)=\mu(0)\left( \frac{b_0+t}{b_0} \right)^{b_1/2}. \ee 
Choose
$\mu(0)$ so that \eqref{e30} holds, and fix it.  Then
inequality \eqref{e31} holds if 
\be\label{e40} 2c_0\mu^{-p}(0)\leq
\frac{b_1}{b_0+t}\frac{(b_0+t)^{\frac{b_1p}{2}}}{b_0^{\frac{b_1p}{2}}},\quad
\forall t\geq 0. \ee 
Choose $b_1$ so that \be\label{e41}
b_1p>2,\qquad p>0. \ee 
Then \eqref{e40} holds if and only if it holds for $t=0$, that is:
\be\label{e42} 2c_0\mu^{-p}(0)\leq
\frac{b_1}{b_0}. \ee Inequality \eqref{e42} is
satisfied if either $b_1$ is chosen sufficiently large for any fixed
$b_0$, or $b_0$ is chosen sufficiently small for any fixed
$b_1>2p^{-1}$ (see \eqref{e41}). In either case one concludes that
the zero solution to equation \eqref{e1} is globally asymptotically
stable.\\
\thmref{thm1} is proved. \hfill $\square$

\section{Additional results. Examples}
{\bf Example 1.} Consider two equations: \be\label{e43} \dot{u}(t)=Au(t),
\ee
\be\label{e44} \dot{v}(t)=Av(t)+B(t)v(t),\qquad t\geq 0, \ee where $A$
and $B(t)$ are bounded linear operators in $H$, $A$ {\it does not depend
on $t$}, and \be\label{e45} \int_0^\infty\|B(t)\|dt<\infty. \ee We
assume that all the solutions to \eqref{e43} are bounded. Then
by the Banach-Steinhaus theorem the following inequality holds:
\be\label{e46} \sup_{t\geq 0}\|e^{tA}\|\leq c<\infty. \ee This
implies Lyapunov's stability of the zero solution to \eqref{e43},
and the inclusion $\sigma(A)\subset \sqcap:=\{z\ : \ \text{Re} z\leq
0\}$, which implies Re$(Au,u)\leq 0$ $\forall u\in H$. A well-known
result is (see, e.g., \cite{B}):

{\it If \eqref{e45} and \eqref{e46} hold
then the zero solution to \eqref{e44} is Lyapunov stable.}

The usual
proof (see \cite{B}, where $H=\R^n$) is based on the Gronwall
inequality. We give a new simple proof based on \lemref{lem1}. Let
$g(t):=\|v(t)\|$. Multiply \eqref{e44} by $u$, take the real part and
use the inequality  Re$(Av,v)\leq 0$ to get: $g\dot{g}\leq 
\|B(t)\|g^2(t)$, $t\geq 0$. Using the inequalities  $g(t)\geq 0$ and 
 \eqref{e45}, one obtains
\be\label{e47} \dot{g}(t)\leq \|B(t)\|g(t),\quad g(t)\leq
g(0)e^{\int_0^\infty\|B(s)\|ds}:=c_1g(0). \ee Therefore, the zero
solution to \eqref{e44} is Lyapunov stable. Moreover, since
$|\dot{g}(t)|\in L^1(\R_+)$, it follows that there exists the finite 
limit:  $\lim_{t\to \infty}\|v(t)\|:=V.$

{\bf Example 2.} Consider a theorem
of N. Levinson in $\R^n$ (see \cite{L} and \cite{D}, pp. 159-164):

{\it If \eqref{e45} and \eqref{e46} hold, then  for every solution $v$ to 
\eqref{e44} one can find a
solution $u$ to \eqref{e43} such that \be\label{e48}\lim_{t\to
\infty}\|u(t)-v(t)\|=0 .\ee}

We give  a new short proof of a
generalization of this theorem to an infinite-dimensional Hilbert
space $H$. If \eqref{e45} and \eqref{e46} hold, then, as we have
proved in Example 1, $\sup_{t\geq 0}\|v(t)\|<\infty$, $\sup_{t\geq
0}\|u(t)\|<\infty$. If $u(0)=u_0$, then $u(t)=e^{tA}u_0$ solves
\eqref{e43}.  Let $v(t)$ solve the equation
\be\label{e49} v(t)=e^{tA}u_0-\int_t^\infty e^{(t-s)A}B(s)v(s)ds.\ee
A simple calculation shows that $v(t)$ solves \eqref{e44} and 
\be\label{e50} \|v(t)-u(t)\|\leq
\int_t^\infty\|e^{(t-s)A}\|\|B(s)\|\|v(s)\|ds\leq
C\int_t^\infty\|B(s)\|ds\to 0,\ t\to \infty, \ee where
$$C=\sup_{t\geq 0}\|e^{tA}\|\sup_{t\geq 0}\|v(t)\|<\infty.$$ The
generalization of Levinson's theorem for $H$ is proved. \hfill
$\square$\\

Equation \eqref{e49} is uniquely solvable in $H$ by iterations for all
sufficiently large $t$ because for such $t$ the norm of the integral
operator in \eqref{e49} is less than one. The unique solution to
\eqref{e49} for sufficiently large $t$ defines uniquely the solution $v$
to \eqref{e44} which satisfies \eqref{e48}.

{\bf Remark 2.}
{\it Our methods are applicable to the equation \eqref{e1} with a force
term: $ \dot{u}=A(t)u+F(t,u)+f(t).$ 
}

\end{document}